\documentclass{amsart}
\usepackage{amsfonts}
\usepackage{amsmath,amssymb}
\usepackage{amsthm}
\usepackage{amscd}
\usepackage{graphics}
\usepackage{graphicx}
\theoremstyle{definition}{
\newtheorem{Def}{{\rm Definition}}
\newtheorem{Ex}{{\rm Example}}
\newtheorem{Rem}{{\rm Remark}}

}
\theoremstyle{plain}
{
\newtheorem{Cor}{Corollary}
\newtheorem{Prop}{Proposition}
\newtheorem{Thm}{Theorem}
\newtheorem{MainThm}{Main Theorem}
\newtheorem{MainCor}{Main Corollary}

}

\begin{document}
\title[Special generic maps on simply-connected manifolds into ${\mathbb{R}}^5$]{Restrictions on special generic maps into ${\mathbb{R}}^5$ on $6$-dimensional or higher dimensional closed and simply-connected manifolds}
\author{Naoki Kitazawa}
\address{Institute of Mathematics for Industry, Kyushu University, 744 Motooka, Nishi-ku Fukuoka 819-0395, Japan\\
 TEL (Office): +81-92-802-4402 \\
 FAX (Office): +81-92-802-4405 \\
}
\email{n-kitazawa@imi.kyushu-u.ac.jp}
\urladdr{https://naokikitazawa.github.io/NaokiKitazawa.html}

\begin{abstract}
The class of {\it special generic} maps is a natural class of smooth maps containing
%naturaｌ class of smooth maps containing
Morse functions on spheres with exactly two singular points and canonical projections of unit spheres. 
We find new restrictions on such maps on $6$-dimensional or higher dimensional closed and simply-connected manifolds into ${\mathbb{R}}^5$. 

Spheres which are not diffeomorphic to unit spheres do not admit such maps whose codimensions are negative in considerable cases. They restrict the homeomorphism and the diffeomorphism types of the manifolds in general. On the other hands, some elementary manifolds admit special generic maps into suitable Euclidean spaces: manifolds represented as connected sums of products of unit spheres are of such examples. This motivates us to study the (non-)existence of special generic maps on elementary manifolds such as projective spaces and some closed and simply-connected manifolds. For example, new explicit investigations of cohomology rings are keys in our new study.

%  For example, we investigate algebraic topological restrictions on {\it Reeb spaces} of the fold maps. %They are essential tools in studying manifolds by using generic maps, and the source manfolds. We also %show flexibility of homology groups of the Reeb spaces and the source manifolds, for example.    

% Abstract text, usually no more than 200 words.
% Avoid bibliographic references (\cite) and complicated mathematics.
% Please do not use custom macros here, as this abstract has to 
% be able to stand alone.  You may use standard tex/latex/AMS macros.
\end{abstract}

% Leave these items like this, and the journal will fill them in.
% \received{Month Day, Year}   % receive date (for example: October 11, 1999)
% \revised{Month Day, Year}    % date of revision; omit, if no revision;
%                             % if multiple revisions, separate by commas
% \published{Month Day, Year}  % publish date\submitted{Bill Murray}      % Name of Journal's Editor, 
% who handled Article 
% \volumeyear{2014} % Volume Year
% \volumenumber{16} % Volume Number 
% \issuenumber{2}   % Issue Number
% \startpage{1}     % PageNumber of first page
% \articlenumber{1} % Sequence number of article within issue
% If copyright is retained by author, comment this out:
% \owner{International Press}

\maketitle
\section{Introduction.}
\label{sec:1}
  
So-called {\it Morse} functions on spheres with exactly two singular points and canonical projections of unit spheres are regarded as {\it special generic} maps. {\it Special generic} maps are also smooth maps whose codimensions are not positive.

These maps have been shown to restrict the differentiable structures of the spheres and more generally restrict the topologies and the differentiable structures of the manifolds strongly.

On the other hands, special generic maps seem to represent manifolds naturally in some senses. Manifolds represented as connected sums of products of unit spheres and other manifolds seeming to be elementary in geometry admit natural special generic maps into suitable Euclidean spaces.

This nature of special generic maps has attracted us in studies of algebraic topological or differential topological properties theory of smooth manifolds. The existence of special generic maps is a natural and attractive problem. It is also difficult in general due to difficulty in construction of explicit smooth maps and finding meaningful restrictions on the manifolds for example. 

We introduce Main theorems. Before this, we introduce notation on algebras and (co)homology. $\mathbb{Z}$ denotes the ring of all integers and $\mathbb{Q} (\supset \mathbb{Z})$ denotes the ring of all rational numbers. 
Let $(X,X^{\prime})$ be a pair of topological spaces satisfying $X^{\prime} \subset X$ where $X^{\prime}$ may be empty. The homology group of the pair $(X,X^{\prime})$ of topological spaces satisfying $X^{\prime} \subset X$ whose coefficient ring is $A$ is denoted by $H_{\ast}(X,X^{\prime};A)$ and the cohomology group of the pair $(X,X^{\prime})$ whose  coefficient ring is $A$ is denoted by $H^{\ast}(X,X^{\prime};A)$. If $A$ is isomorphic to $\mathbb{Z}$ (resp. $\mathbb{Q}$), then the homology group and the cohomology group are called the {\it integral} (resp. {\it rational}) {\it homology group} and the {\it integral} (resp. {\it rational}) {\it cohomology group}, respectively. The $k$-th homology group (cohomology group) is denoted by $H_k(X,X^{\prime};A)$ (resp. $H^k(X,X^{\prime};A)$). If $X^{\prime}$ is empty, then we omit ",$X^{\prime}$" in the notation here and the homology group (cohomology group) of the pair $(X,X^{\prime})$ is also called the homology group (resp. cohomology group) of $X$.
(Co)homology classes of $(X,X^{\prime})$ (or $X$) are elements of the (resp. co)homology groups.

We add "integral" and "rational" as before for these (co)homology groups and (co)classes according to the coefficient rings.

For a topological space $X$, the $k$-th homotopy group of $X$ is denoted by ${\pi}_k(X)$.

Let $(X,X^{\prime})$ and $(Y,Y^{\prime})$ be pairs of topological spaces satisfying $X^{\prime} \subset X$ and $Y^{\prime} \subset Y$ where $X^{\prime}$ and $Y^{\prime}$ can be empty. For a continuous map $c:X \rightarrow Y$ satisfying $c(X^{\prime}) \subset Y^{\prime}$, $c_{\ast}:H_{\ast}(X,X^{\prime};A) \rightarrow H_{\ast}(Y,Y^{\prime};A)$ and $c^{\ast}:H^{\ast}(Y,Y^{\prime};A) \rightarrow H^{\ast}(X,X^{\prime};A)$ and $c_{\ast}:{\pi}_k(X) \rightarrow {\pi}_k(Y)$ denote canonically induced homomorphisms where for homology groups and homotopy groups we abuse the same notation. 

Let $X$ be a topological space. For a sequence $\{a_j\}_{j=1}^l \subset H^{\ast}(X;A)$ of cohomology classes of length $l>0$, ${\cup}_{j=1}^l a_j$ denotes the cup product for general $l$ and $a_1 \cup a_2$ denotes for $l=2$. This makes $H^{\ast}(X;A)$ a graded commutative algebra and this is the cohomology ring of $X$ whose coefficient ring is $A$.
 
For related introductory explanations and some advanced ones on algebraic topology, see \cite{hatcher} as a textbook for example.

We present our Main Theorems and Main Corollary. We leave rigorous and precise expositions on the notation and notions we need later. In short, a special generic map $f:M \rightarrow {\mathbb{R}}^n$ from an $m$-dimensional closed smooth manifold $M$ into the $n$-dimensional Euclidean space ${\mathbb{R}}^n$ satisfying $m>n$ is known to be represented as the composition of a smooth surjection $q_f$ onto an $n$-dimensional compact smooth manifold $W_f$ with a smooth immersion $\bar{f}:W_f \rightarrow {\mathbb{R}}^n$ (Proposition \ref{prop:2} (\ref{prop:2.1})). The {\it singular set} $S(f) \subset M$ of the special generic map $f$ is defined as the set of all points in the manifold $M$ where the ranks of the differentials are smaller than $n$ and a closed smooth submanifold with no boundary of dimension $n-1$ (Proposition \ref{prop:1}).
\begin{MainThm}
\label{mthm:1}
Suppose that a special generic map $f:M \rightarrow {\mathbb{R}}^5$ on an $m$-dimensional closed and simply-connected manifold $M$ exists where $m \geq 6$. Assume also that a family $\{h_j\}_{j=1}^l \subset H_2(M;\mathbb{Z})$ of homology classes of length $l>0$ satisfying the following conditions exists and that $H_2(M;\mathbb{Z})$ and $H_2(W_f;\mathbb{Z})$ are of rank $l$.
\begin{itemize}

\item Each $h_j$ is not divisible by any integer greater than $1$.
\item Each $h_j$ is of infinite order.
\item The $l>0$ homology classes in the family are mutually independent.
\item For each $h_j$, we can take a smooth embedding $e_j:S^2 \rightarrow M$ of the $2$-dimensional unit sphere $S^2$ satisfying the following properties.
\begin{itemize}
\item Each $h_j$ is represented by the submanifold $e_j(S^2)$ for a suitable orientation of $S^2$.
\item A homology class ${h_j}^{\prime} \in H_2(W_f;\mathbb{Z})$ of infinite order can be taken for each $1 \leq j \leq l$.
\item Each ${h_j}^{\prime}$ is realized as the value of the homomorphism ${(q_f \circ e_j)}_{\ast}:H_2(S^2;\mathbb{Z}) \rightarrow H_2(W_f;\mathbb{Z})$ at the fundamental class for the suitable orientation of $S^2$ before.
\item The $l>0$ elements in the family $\{{h_j}^{\prime}\}_{j=1}^l$ are mutually independent.
\end{itemize}
%\item Points in the sphere $S^2$ are mapped by the map $q_f \circ e_j$ into the interior ${\rm Int}\ W_f \subset W_f$.
\end{itemize} 
Then we have a family $\{{h_j}^{\rm c}\}_{j=1}^l \subset H^2(M;\mathbb{Z})$ of cohomology classes of length $l>0$ satisfying the following conditions.
\begin{enumerate}
\item
\label{mthm1:1}
Each ${h_j}^{\rm c}$ is not divisible by any integer greater than $1$.
\item
\label{mthm1:2}
The $l>0$ cohomology classes in the family are mutually independent.
\item
\label{mthm:1.3}
The cup product of ${h_{j_1}}^{\rm c}$ and ${h_{j_2}}^{\rm c}$ satisfies the following properties.
\begin{enumerate}
\item
\label{mthm:1.3.1}
Suppose also that the {\rm singular set} $S(f)$ of $f$ is connected. Then it vanishes for any integers $1 \leq j_1,j_2 \leq l$.
\item
\label{mthm:1.3.2}
Suppose also that each ${h_j}^{\prime} \in H_2(W_f;\mathbb{Z})$ is not divisible by $2$. By using the canonical homomorphism defined as the quotient map ${\phi}_{2,\mathbb{Z}}$ from $\mathbb{Z}$ onto the group $\mathbb{Z}/2\mathbb{Z}$ of order $2$, which is also a field, we have ${\phi}_{2,\mathbb{Z}}({h_{j}}^{\rm c}) \in H^2(M;\mathbb{Z}/2\mathbb{Z})$ for any $1 \leq j \leq l$. The square of ${\phi}_{2,\mathbb{Z}}({h_{j}}^{\rm c}) \in H^2(M;\mathbb{Z}/2\mathbb{Z})$ always vanishes.
\end{enumerate}
\end{enumerate}

\end{MainThm}
\begin{MainCor}
\label{mcor:1}
Let $M$ be the $3$-dimensional complex projective space. If $M$ admits a special generic map $f:M \rightarrow {\mathbb{R}}^5$, then $M$ does not have the family $\{h_j\}_{j=1}^l \subset H_2(M;\mathbb{Z})$ of the homology classes as in Main Theorem \ref{mthm:1}.
\end{MainCor}
\begin{MainCor}
	\label{mcor:2}
	Let $m \geq 7$ be an integer.
	Let $M$ be an $m$-dimensional closed and simply-connected manifold. If there exists a 2nd integral cohomology class whose square is not divisible by $2$, then $M$ admits no special generic maps into $ {\mathbb{R}}^5$.
\end{MainCor}
\begin{MainThm}
\label{mthm:2}
Suppose that a special generic map $f:M \rightarrow {\mathbb{R}}^5$ on an $m$-dimensional closed and simply-connected manifold $M$ exists where $m \geq 6$. Assume also that a homology class $h_0 \in H_2(M;\mathbb{Z})$ of infinite order exists. We have the following two.
\begin{enumerate}
\item $h_0$ is represented by a closed smooth submanifold $S$ and we may regard this submanifold being diffeomorphic to $S^2$ by virtue of Hurewicz theorem.
\item We can not choose an embedding $e:S^2 \rightarrow M$ giving some submanifold $S=e(S^2)$ diffeomorphic to $S^2$ satisfying both the condition before and all the following properties.
\begin{enumerate}
\item The image of the composition $q_f \circ e$ is a subset of ${\rm Int}\ W_f$.
\item A homology class of $H_2(W_f;\mathbb{Z})$ of a finite order is realized as the value of the homomorphism $(q_f \circ e)_{\ast}:H_2(S^2;\mathbb{Z}) \rightarrow H_2(W_f;\mathbb{Z})$ at the fundamental class of $S^2$ of the domain for any orientation of it.
\end{enumerate}
\end{enumerate}
\end{MainThm}
Although $m \geq 6$ is assumed, this is essentially a theorem for $m=6$. Remark \ref{rem:4} explains about this. This is also a more general remark related to the present paper.

Main Theorems \ref{mthm:1} and \ref{mthm:2} with additional arguments yield our last Main Theorem.

\begin{MainThm}
\label{mthm:3}
Assume that a special generic map $f:M \rightarrow {\mathbb{R}}^5$ on an $n$-dimensional closed and simply-connected manifold $M$ whose singular set $S(f)$ is connected exists where $m \geq 7$. Then the cup product $c_1 \cup c_2$ of any pair of cohomology classes $c_1,c_2 \in H^2(M;\mathbb{Z})$ always vanishes.
\end{MainThm}

It has been announced in \cite{kitazawa1} by the author that real projective spaces do not admit special generic maps whose codimensions $m-n$
are negative (Theorem \ref{thm:3} and Corollary \ref{cor:1}). In the case where the codimension is $m-n=0$ the result depends on the dimension of the projective space. 
Theorem \ref{thm:3} also implies that the existence of a special generic map $f:M \rightarrow {\mathbb{R}}^n$ on an $m$-dimensional closed and simply-connected manifold $M$ into the $n$-dimensional Euclidean space ${\mathbb{R}}^n$ yields the vanishing of the cup product $c_1 \cup c_2$ of any pair of two cohomology classes $c_1,c_2 \in H^2(M;\mathbb{Z})$ for $m \geq 6$ and $n=1,2,3,4$ and this implies that Main theorem \ref{mthm:3} is a version for $n=5$ of this. Special generic maps on closed and simply-connected manifolds were studied first in \cite{saeki}, followed by \cite{nishioka, wrazidlo} and works \cite{kitazawa2, kitazawa3} of the author for example. Theorems \ref{thm:1} and \ref{thm:2} present some of them. These studies are essentially on such maps on spheres whose dimensions are arbitrary and closed and simply-connected manifolds whose dimensions are at most $5$. Cases where the manifolds are of dimensions greater than $5$ are new except several facts we immediately have from existing studies. As a kind of appendices, Remark \ref{rem:5} refers to classifications of $6$-dimensional closed and simply-connected manifolds where \cite{barden} is a $5$-dimensional version for example. Closed simply-connected manifolds whose dimensions are at least $5$ have been classified via algebraic and abstract systems in the midst of the 20th century.

The present paper is organized as follows. The next section is for preliminary. We also review the definition of a special generic map (Definition \ref{def:1}). After the next section, we review and introduce some existing studies of special generic maps including ones before. In the fourth section, we prove our Main Theorems.

%{\it Fold} maps are, in short, higher dimensional versions of Morse functions and important tools in studying geometric properties of manifolds in the branch of the singularity theory of differentiable maps and application to geometry of manifolds as Morse functions are in so-called Morse theory: Morse theory is in a sense regarded as specific classical theory of the singularity theory of differentiable maps. $7$-dimensional manifolds form an interesting class of manifolds. Especially, Milnor's discovery of $7$-dimensional {\it exotic} homotopy spheres, which are homotopy spheres having differentiable structures different from those of standard spheres \cite{milnor}, has made this class more attractive and this class is still attractive. Recent attractive studies on 
\section{Preliminary.}
${\mathbb{R}}^k$ denotes the $k$-dimensional Euclidean space.

$\mathbb{R}^1$ is denoted by $\mathbb{R}$ in usual situations. $\mathbb{R} \supset \mathbb{Q}$ of course. However we do not need these in the present paper essentially.
 
We assume that this space is naturally a $k$-dimensional smooth manifold endowed with the standard Euclidean metric. For any point $x \in {\mathbb{R}}^k$, $||x||$ denotes the distance between $x$ and the origin $0 \in {\mathbb{R}}^k$ where the metric is considered. For an arbitrary integer $k \geq 1$, $D^k:=\{x \mid ||x|| \leq 1\} \subset {\mathbb{R}}^k$ denotes the $k$-dimensional unit disk and this is a $k$-dimensional compact, connected and smooth closed submanifold. For an arbitrary integer $k \geq 0$, $S^k:=\{x \mid ||x||=1\} \subset {\mathbb{R}}^{k+1}$ denotes the $k$-dimensional unit sphere. This is a $k$-dimensional closed smooth submanifold with no boundary. This is also a two-point discrete set for $k=0$ and for $k \geq 1$ it is connected.

A smooth manifold is well-known to have the structure of a canonical PL manifold. We regard smooth manifolds as such PL manifolds in several scenes in the present paper.
It is also well-known that topological manifolds are regarded as CW complexes.

A {\it singular} point $p \in X$ of a smooth map $c:X \rightarrow Y$ is a point where the rank of the differential ${dc}_p$ is smaller than $\min\{\dim X,\dim Y\}$. The {\it singular set} $S(c)$ of $c$ is the set of all singular points of $c$.

A {\it diffeomorphism} is a smooth map between two manifolds which is also a homeomorphism with no singular points. A diffeomorphism from a manifold $X$ into the same manifold is also said to be a {\it diffeomorphism} on $X$. The diffeomorphism group of a smooth manifold $X$ is the group of all diffeomorphisms on $X$, topologized with the so-called {\it Whitney $C^{\infty}$ topology}. For this topology, see \cite{golubitskyguillemin} for example. 

\cite{golubitskyguillemin} is also on fundamental or advanced theory of differentiable or smooth maps and related fundamental or advanced theory on differential topology of manifolds. Methods and theory presented there are also applied to our proof of our Main Theorem and others explicitly or implicitly. {\it Generic} smooth immersions and embeddings, which we can also consider in the other categories such as the PL category naturally, appear for example.

On the family of all smooth manifolds, we can define equivalence relations by the following ways. 
\begin{itemize}
\item Two manifolds are equivalent if they are mutually homeomorphic. For this, we call an equivalence class a {\it homeomorphism type}.
\item Two manifolds are equivalent if there exists a diffeomorphism from a manifold into the other manifold. For this, we call an equivalence class a {\it diffeomorphism type}. We also say that they are mutually {\it diffeomorphic}.  
\end{itemize}

If we consider orientations on the manifolds, then we add "{\it oriented}". We do not consider {\it oriented homeomorphism types} in the present paper and we only consider {\it oriented diffeomorphism types} for oriented smooth manifolds.  
A {\it homotopy sphere} is a smooth manifold whose homeomorphism type is same as that of a unit sphere.
A {\it standard} sphere is a homotopy sphere whose diffeomorphism type is same as that of a unit sphere.
An {\it exotic} sphere is a homotopy sphere whose diffeomorphism type is not same as that of any unit sphere.
%The {\it singular value set} of $c$ is the image $c(S(c))$ of the singular set. The {\it regular value set} of $c$ is the complementary set $Y-c(S(c))$. A {\it singular {\rm (}regular{\rm )} value} is a point in the singular (resp. regular) value set of this map. 

\begin{Def}
\label{def:1}
Let $m \geq n \geq 1$ be integers. A smooth map from an $m$-dimensional closed smooth manifold into ${\mathbb{R}}^n$ is said to be {\it special generic} if at each singular point $p$, it is locally represented by the form
$$(x_1, \cdots, x_m) \mapsto (x_1,\cdots,x_{n-1},\sum_{k=n}^{m}{x_k}^2)$$
for suitable coordinates. 
\end{Def}
\begin{Rem}
\label{rem:1}
We can replace ${\mathbb{R}}^n$ by a general $n$-dimensional smooth manifold with no boundary satisfying some conditions throughout the present paper. However, we concentrate on cases where the manifolds of the targets are ${\mathbb{R}}^n$.
\end{Rem}
\begin{Prop}
\label{prop:1}
The singular set of a special generic map in Definition {\rm \ref{def:1}} is an {\rm (}$n-1${\rm )}-dimensional smooth closed submanifold with no boundary and the restriction of the map there is a smooth immersion.
\end{Prop}

\section{Some existing studies on special generic maps.}

A bundle whose fiber is a smooth manifold and whose structure group is a subgroup of the diffeomorphism group is said to be a {\it smooth} bundle. A smooth bundle whose fiber is an Euclidean space, a unit sphere, or a unit disk and whose structure group consists of (natural) linear transformations is said to be {\it linear}.

We implicitly apply fundamental or advanced arguments on bundles in the present paper. For them, see \cite{milnorstasheff, steenrod} for example.  

\begin{Prop}[E. g. \cite{saeki}]
\label{prop:2}
Let $m>n \geq 1$ be integers.
\begin{enumerate}
\item
\label{prop:2.1}
For a special generic map $f:M \rightarrow {\mathbb{R}}^n$ on an $m$-dimensional closed and connected manifold $M$, the following properties hold.
\begin{enumerate}
\item There exists an $n$-dimensional compact manifold $W_f$ smoothly immersed into ${\mathbb{R}}^n$ via $\bar{f}:W_f \rightarrow {\mathbb{R}}^n$.
\item There exists a smooth surjection $q_f:M \rightarrow W_f$.
\item $q_f$ maps the singular set $S(f)$ onto the boundary $\partial W_f \subset W_f$ as a diffeomorphism.
\item We can choose a small collar neighborhood $N(\partial W_f) \subset W_f$ and the composition of the map $q_f {\mid}_{{q_f}^{-1}(N(\partial W_f))}$ onto $N(\partial W_f)$ with the canonical projection to $\partial W_f$ gives a linear bundle whose fiber is diffeomorphic to the unit disk $D^{m-n+1}$.
\item The restriction of $q_f$ to the preimage of $W_f-{\rm Int}\ N(\partial W_f)$ gives a smooth bundle over $W_f-{\rm Int}\ N(\partial W_f)$ whose fiber is diffeomorphic to the unit sphere $S^{m-n}$.
\end{enumerate} 
\item
\label{prop:2.2}
For an $n$-dimensional compact, connected and smooth manifold $W_N$ and a smooth immersion ${\bar{f}}_N:W_N \rightarrow {\mathbb{R}}^n$, a special generic map $f_0:M_0 \rightarrow {\mathbb{R}}^n$ on some $m$-dimensional closed and connected manifold $M_0$ exists and the following properties hold where we abuse similar notation as the previous proposition {\rm (\ref{prop:2.1})}.
\begin{enumerate}
\item There exists a diffeomorphism $\phi:W_{f_0} \rightarrow W_N$ such that $\bar{f_0}={\bar{f}}_N \circ \phi$.
\item We can choose a small collar neighborhood $N(\partial W_{f_0}) \subset W_{f_0}$ and the composition of the map $q_{f_0} {\mid}_{{q_{f_0}}^{-1}(N(\partial W_{f_0}))}$
onto $N(\partial W_{f_0})$ with the canonical projection to $\partial W_{f_0}$ gives a trivial linear bundle whose fiber is diffeomorphic to the unit disk $D^{m-n+1}$.
\item The restriction of $q_{f_0}$ to the preimage of $W_{f_0}-{\rm Int}\ N(\partial W_{f_0})$
gives a trivial smooth bundle over $W_{f_0}-{\rm Int}\ N(\partial W_{f_0})$
whose fiber is diffeomorphic to the unit sphere $S^{m-n}$.
\end{enumerate} 
\end{enumerate}
\end{Prop}

For a finite set $X$, $|X|$ denotes the size of $X$.
\begin{Ex}
\label{ex:1}
Let $m \geq n \geq 2$ be integers. Let $\{S^{k_j} \times S^{m-k_j}\}_{j \in J}$ be a family of finitely many products of two unit spheres where $k_j$ is an integer satisfying $1 \leq k_j \leq n-1$. We consider a connected sum of these $J$ manifolds in the smooth category and set $M$ as an $m$-dimensional closed and connected manifold diffeomorphic to this. 
\begin{enumerate}
\item \label{ex:1.1}
Let $m>n$. Then $M$ admits a special generic map $f$ satisfying the following properties.
\begin{enumerate}
\item $f {\mid}_{S(f)}$ is an embedding.
\item $W_f$ is regarded as the image $f(M)$ in a canonical way. The image $W_f$ is diffeomorphic to a manifold obtained as a boundary connected sum of $|J|$ manifolds in the family $\{S^{k_j} \times D^{n-k_j}\}_{j \in J}$. Of course the boundary connected sum is considered in the smooth category.
\item We can choose a small collar neighborhood $N(\partial W_f) \subset W_f$ and the composition of the map $q_f {\mid}_{{q_f}^{-1}(N(\partial W_f))}$ onto $N(\partial W_f)$ with the canonical projection to $\partial W_f$ gives a trivial linear bundle whose fiber is diffeomorphic to the unit disk $D^{m-n+1}$.
\item The restriction of $q_f$ to the preimage of $W_f-{\rm Int}\ N(\partial W_f)$ gives a smooth bundle over $W_f-{\rm Int}\ N(\partial W_f)$ whose fiber is diffeomorphic to $S^{m-n}$.
\end{enumerate}
This is one of simplest examples for Proposition \ref{prop:2}.
\item \label{ex:1.2}
Let $m=n$. Then $M$ admits a special generic map $f$ satisfying the following properties.
\begin{enumerate}
\item $f {\mid}_{S(f)}$ is an embedding.
\item The image $f(M)$ is diffeomorphic to a manifold obtained as a boundary connected sum of $|J|$ manifolds in the family $\{S^{k_j} \times D^{n-k_j}\}_{j \in J}$. Of course the boundary connected sum is considered in the smooth category.
\item We can choose a small collar neighborhood $N(\partial f(M)) \subset f(M)$ and the composition of the map $f {\mid}_{{f}^{-1}(N(\partial f(M)))}$ onto $N(\partial f(M))$ with the canonical projection to $\partial f(M)$ gives a trivial linear bundle whose fiber is diffeomorphic to the unit disk $D^{m-n+1}=D^1$.
\item The restriction of $f$ to the preimage of $f(M)-{\rm Int}\ N(\partial f(M))$ gives a smooth bundle over $f(M)-{\rm Int}\ N(\partial f(M))$ whose fiber is diffeomorphic to $S^{m-n}=S^0$.
\end{enumerate}
\end{enumerate}
\end{Ex}
\begin{Prop}[E. g. \cite{saeki}]
\label{prop:3}
Let $m>n \geq 1$ be integers. 
\begin{enumerate}
\item For a special generic map $f:M \rightarrow {\mathbb{R}}^n$ on an $m$-dimensional closed and connected manifold $M$, there exists an {\rm (}$m+1${\rm )}-dimensional compact, connected and {\rm (}{\rm PL}{\rm )} manifold $W$ whose boundary is $M$ and which collapses to $W_f$, regarded as a subpolyhedron of $W$ in a suitable way. If we consider PL manifolds here, then we consider boundaries in the PL category. Furthermore, for the canonical inclusion $i_M:M \rightarrow W$ and a continuous {\rm (}resp. PL{\rm )} map $r:W \rightarrow W_f$ giving a collapsing to $W_f$, we have $q_f=r \circ i_M$. If $m-n=1,2,3$ in addition, then we can take $W$ as a smooth manifold and $r$ as a smooth map. 
\item In Proposition {\rm \ref{prop:2} (\ref{prop:2.2})}, there exists a suitable construction of the map $f_0:M_0 \rightarrow {\mathbb{R}}^n$ enabling us to take an {\rm (}$m+1${\rm )}-dimensional compact, connected and smooth manifold $W_0$ whose boundary is diffeomorphic to $M_0$ and which collapses to $W_{f_0}$, regarded as a subpolyhedron of $W_0$ in a suitable way. Furthermore, for the canonical inclusion $i_{M,0}:M_0 \rightarrow W_0$ and a PL map $r_0:W_0 \rightarrow W_{f_0}$ giving a collapsing to $W_{f_0}$, we have $q_{f_0}=r_0 \circ i_{M,0}$. Furthermore, $r_0$ can be taken as a smooth map. 
\end{enumerate}
\end{Prop}

For more general propositions of this type, see \cite{saekisuzuoka} and see papers \cite{kitazawa0.1,kitazawa0.2,kitazawa0.3} by the author.

The following theorem shows existing results on the differential structures or the (oriented) diffeomorphism types of homotopy spheres admitting special generic maps. 
A closed manifold admits a special generic function or a special generic map for $n=1$ in Definition \ref{def:1} if and only if it is a homotopy sphere which is not a $4$-dimensional exotic sphere. This is due to so-called Reeb's theorem and theory of Morse functions. For classical important theory on Morse functions, see \cite{milnor2,milnor3} for example. For (oriented) diffeomorphism types of homotopy spheres, see \cite{milnor} as a pioneering study and see also \cite{eellskuiper,kervairemilnor} for example.

\begin{Thm}[\cite{calabi, saeki, saeki2, wrazidlo}]
\label{thm:1}
Let $m>n \geq 1$ be integers.
\begin{enumerate}
\item A special generic map $f:M \rightarrow {\mathbb{R}}^n$ on an $m$-dimensional closed and connected manifold $M$ is a homotopy sphere if and only if $W_f$ is contractible. 
\item \label{thm:1.2} Let $m>1$. Every $m$-dimensional homotopy sphere which is not a $4$-dimensional exotic sphere admits a special generic map into the plane such that restriction to the singular set is an embedding and that the singular value set is a smoothly embedded circle. $4$-dimensional exotic spheres do not admit special generic maps into ${\mathbb{R}}^n$ for $n=1,2$ and note that such manifolds are still undiscovered.
\item \label{thm:1.3} Let $m>3$. $m$-dimensional exotic homotopy spheres do not admit special generic maps into ${\mathbb{R}}^{m-3}$, ${\mathbb{R}}^{m-2}$ and ${\mathbb{R}}^{m-1}$.
\item \label{thm:1.4} Let $m=7$. $m$-dimensional oriented homotopy spheres of $14$ types of all $28$ oriented diffeomorphism types do not admit special generic maps into ${\mathbb{R}}^3$.
\end{enumerate}
\end{Thm}
In addition, every $m$-dimensional homotopy sphere which is not a $4$-dimensional exotic sphere admits a special generic map into ${\mathbb{R}}^m$ according to \cite{eliashberg}.
\begin{Thm}
\label{thm:2}
Let $m \geq n \geq 1$ be integers{\rm :} connected sums are considered in the smooth category.
\begin{enumerate}
\item {\rm (\cite{saeki})}
\label{thm:2.1}
Let $m \geq 2$. An $m$-dimensional closed and connected manifold admits a special generic map into ${\mathbb{R}}^2$ if and only if it is a homotopy sphere which is not a $4$-dimensional exotic sphere or represented as a connected sum of the total spaces of smooth bundles over $S^1$ whose fibers are diffeomorphic to homotopy spheres which are not $4$-dimensional exotic spheres.
\item {\rm (\cite{saeki})}
\label{thm:2.2}
Let $m=4,5$. An $m$-dimensional closed and simply-connected manifold admits a special generic map into ${\mathbb{R}}^3$ if and only if it is a standard sphere or represented as a connected sum of the total spaces of smooth bundles over $S^2$ whose fibers are diffeomorphic to $S^3$.
\item {\rm (\cite{nishioka})}
\label{thm:2.3}
Let $m=5$. An $m$-dimensional closed and simply-connected manifold admits a special generic map into ${\mathbb{R}}^n$ for $n=1,2$ if and only if it is a standard sphere. An $m$-dimensional closed and simply-connected manifold admits a special generic map into ${\mathbb{R}}^n$ for $n=3, 4$ if and only if it is a standard sphere or represented as a connected sum of the total spaces of smooth bundles over $S^2$ whose fibers are diffeomorphic to $S^3$.
\end{enumerate}
\end{Thm}
In addition, \cite{nishioka} finds a necessary and sufficient condition for a $5$-dimensional closed and simply-connected manifold to admit a special generic map into ${\mathbb{R}}^5$ by applying \cite{barden,eliashberg}. \cite{barden} presents complete classifications of $5$-dimensional closed and simply-connected manifolds in the topology, PL, and smooth categories. These categories are equivalent for $5$-dimensional closed and simply-connected manifolds. \cite{nishioka} finds a necessary and sufficient condition for a $5$-dimensional closed and simply-connected manifold to admit a special generic map into ${\mathbb{R}}^4$ by applying this theory and showing that the 2nd integral homology group of a manifold admitting such a map is free.
\begin{Thm}[\cite{kitazawa1}]
\label{thm:3}
Let $M$ be a closed and connected manifold of dimension $m>1$. Let $n<m$ and $l$ be positive integers.
Let $A$ be a commutative ring. Let there exist a sequence $\{a_j\}_{j=1}^l \subset H^{\ast}(M;A)$ satisfying the following three. 
\begin{itemize}
\item The cup product ${\cup}_{j=1}^l a_j$ does not vanish.
\item The degree of each $a_j$ of the previous classes is smaller than or equal to $m-n$.
\item The sum of the degrees for the previous $l>0$ elements in $\{a_j\}_{j=1}^l$ is greater than or equal to $n$.
\end{itemize}
Then $M$ does not admit special generic maps into ${\mathbb{R}}^n$.
\end{Thm}
We review a proof. We omit erigorous xpositions on definitions of a {\it handle} and its {\it index} for PL manifolds and general polyhedra.
\begin{proof}
Suppose that $M$ admits a special generic map into ${\mathbb{R}}^n$. We can take an ($m+1$)-dimensional compact and connected PL manifold $W$ as in Proposition \ref{prop:3}. 
$W$ is obtained by attaching handle to $M \times \{0\} \subset M \times [-1,0]$ whose indices are greater than $m+1-{\dim W_f}=m-n+1$.
We can take $b_j \in H^{\ast}(W;A)$ uniquely so that $a_j={i_M}^{\ast}(b_j)$ in Proposition \ref{prop:3}. The cup product ${\cup}_{j=1}^k a_j$ vanishes since $W$ collapses to $W_f$, collapsing to an ($n-1$)-dimensional
polyhedron. This contradicts the assumption. This completes the proof.
\end{proof}
\begin{Cor}[\cite{kitazawa1}]
\label{cor:1}
Let $m \geq n \geq 1$ be integers.
\begin{enumerate}
\item Let $m>n$ in addition. The $m$-dimensional real projective space does not admit special generic maps into ${\mathbb{R}}^n$. 
\item Let $n<m-1$ and $m=2k$ for some positive integer $k>1$ in addition. The $m$-dimensional complex projective space, which is also a $2k$-dimensional closed, simply-connected and smooth manifold, does not admit special generic maps into ${\mathbb{R}}^n$.
\end{enumerate}
\end{Cor}
Note that in Corollary \ref{cor:1} in the case $m=n$, \cite{eliashberg} presents a useful tool. 
This theory says a necessary and sufficient condition for a closed and orientable smooth manifold to admit a special generic map into the Euclidean space of the same dimension is that the Whitney sum of the tangent bundle and the trivial real line bundle is trivial. 
We argue the existence of special generic maps whose codimensions are $0$ on real projective spaces only here.
According to this, if $m=n=1,3,7$, then the $m$-dimensional real projective space admits a special generic map into ${\mathbb{R}}^n$ (for $m=n=7$ see also Remark \ref{rem:2}). According to \cite{kikuchisaeki}, if the $m$-dimensional real projective space or more generally, a closed manifold whose Euler number is odd, admits a special generic map into ${\mathbb{R}}^n={\mathbb{R}}^m$, then $m=n=1,3,7$. If $m=n$ is even, then the $m$-dimensional real projective space does not admit special generic maps into ${\mathbb{R}}^n$.

\begin{Rem}
\label{rem:2}
According to a slide \cite{wrazidlo2} of a related talk in a conference, Corollary \ref{cor:1} had been first shown for the $7$-dimensional real projective space by using known facts and theory on $7$-dimensional homotopy spheres and special generic maps in Theorem \ref{thm:1} before \cite{kitazawa1} was announced. After the announcement of \cite{kitazawa1}, \cite{wrazidlo3} announced a proof investigating restrictions on the torsion groups of the integral homology groups for {\it rational homology spheres}: a {\it rational homology sphere} $M$ is a closed and smooth manifold whose rational homology group $H_{\ast}(M;\mathbb{Q})$ is isomorphic to that of a sphere.
\end{Rem}

For other related studies on special generic maps, see \cite{burletderham, furuyaporto,  saekisakuma, saekisakuma2, sakuma} and works \cite{kitazawa2, kitazawa3} of the author for example. \cite{kitazawa4} is also a related preprint. We do not need to understand notions and arguments of the preprint well.

\section{Proofs of our Main Theorems.}
The {\it fundamental class} of a compact, connected and oriented smooth manifold $Y$ is the canonically and uniquely defined ($\dim Y$)-th homology class, which is also a generator of the group $H_{\dim Y}(Y,\partial Y;\mathbb{Z})$, isomorphic to $\mathbb{Z}$. 
Let $i_{Y,X}:Y \rightarrow X$ be a smooth immersion satisfying $i_{Y,X}(\partial Y) \subset \partial X$ and $i_{Y,X}({\rm Int}\ Y) \subset {\rm Int}\ X$. In other words, $Y$ is {\it properly} immersed or embedded into $X$.
If for a homology class $h$ of $X$, the value of the homomorphism ${i_{Y,X}}_{\ast}$ induced by the smooth immersion or embedding $i_{Y,X}:Y \rightarrow X$ at the fundamental class of $Y$ is $h$, then $h$ is said to be {\it represented} by $Y$. We can consider these arguments and notions in the PL category for example.

We present shortly the notions of {\it spin} vector (linear) bundles and {\it spin} manifolds. {\it Spin} bundles are real vector bundles or linear bundles which are orientable and whose {\it 2nd Stiefel-Whitney classes} vanish.
{\it Spin} manifolds are orientable smooth manifolds such that the tangent bundles are spin.
For these notions, consult \cite{milnorstasheff} for example. 

\begin{proof}[A proof of Main Theorem \ref{mthm:1}]
First, according to Proposition \ref{prop:3} and the proof of Theorem \ref{thm:3}, ${\pi}_1(W_f)$ also vanishes. For this see also related arguments in \cite{saekisuzuoka} such as Corollary 4.8 there and articles \cite{kitazawa0.1,kitazawa0.2,kitazawa0.3} by the author.
Another important fact is that homology classes are represented by compact and oriented smooth submanifolds thanks to \cite{thom} in our proof. \\ 

$q_f \circ e_j:S^2 \rightarrow W_f$ can be regarded as a smooth embedding into $W_f$ having finitely many points in $\partial W_f$ by a fundamental argument on singularity theory and differential topology. 
%Let ${e_j}^{\prime}$ denote the embedding.
In addition, we can take a family $\{F_j \subset W_f\}$ of smoothly and properly embedded $3$-dimensional, compact, connected and orientable submanifolds satisfying the following properties from the assumption including the conditions: in other words ${\rm Int}\ F_j$ is embedded into ${\rm Int}\ W_f$ and $\partial F_j$ is embedded into $\partial W_f$ (in a so-called {\it generic} way). 
%Related arguments are also in \cite{kitazawa2} where we do not need essential knowledge about this preprint in our present study.
\begin{itemize}
\item $H_3(W_f;\partial W_f;\mathbb{Z})$ is isomorphic to $H^2(W_f;\mathbb{Z})$, free and of rank $l$ by virtue of the fact that $W_f$ is simply-connected and orientable and Poincar\'e duality theorem for $W_f$. Furthermore, there exists a basis of $H_3(W_f;\partial W_f;\mathbb{Z})$ consisting of exactly $l$ elements satisfying the following two.
\begin{itemize}
	\item Each of these $l$ elements is not divisible by any integer greater than $1$.
	\item Each of these $l$ elements is represented by some $F_j$.
\end{itemize}	
\item ${q_f}^{-1}(F_j)$ is regarded as the ($m-2$)-dimensional closed manifold of the domain of a special generic map into a $3$-dimensional smooth manifold ${F_j}^{\prime}$ with no boundary satisfying ${F_j} \subset F_j^{\prime}$: we consider a special generic map into a general manifold with no boundary here for each $j$. This is essentially thanks to local structures of maps around the singular points and the singular values. 
\item ${q_f}^{-1}(F_{j_0})$ and the union of the images of all embeddings in $\{e_j\}$ intersect in a finite set for each $j_0$. $F_{j_0}$ and the union of the images of all maps in $\{q_f \circ e_j\}$ intersect only in the image of the finite set by the map $q_f$ for each $j_0$.
\item For the family $\{h_j\}_{j=1}^l$, we can define the duals $\{{h_j}^{\ast}\}_{j=1}^l \subset H^2(M;\mathbb{Z})$ satisfying ${h_{j_1}}^{\ast}(h_{j_2})=0$ for $j_1 \neq j_2$ and ${h_j}^{\ast}(h_j)=1$. 
The ($m-2$)-th homology class of $M$ represented by ${q_f}^{-1}(F_j)$ is the so-called {\it Poincar\'e dual} to $a_j {h_j}^{\ast}$ for a suitable non-zero integer $a_j$. 
%\item Furthermore, by virtue of Poincar\'e duality theorem and conditions on $\{{h_j}^{\prime}\} \subset H_2(W_f;\mathbb{Z})$, we can also do so that the cup product of ${h_{j_1}}^{\ast}$ and ${h_{j_2}}^{\ast}$ vanishes for any pair $(j_1,j_2)$ satisfying $j_1 \neq j_2$.
\end{itemize}

We define as ${h_j}^{\rm c}:={h_j}^{\ast}$. We investigate the cup product of ${h_{j_1}}^{\rm c}$ and ${h_{j_2}}^{\rm c}$ to complete the proof. \\
\ \\
\underline{Case 1} $S(f)$ is connected (Main Theorem \ref{mthm:1} (\ref{mthm:1.3.1})). \\
For any pair $(F_{j_1},F_{j_2})$ satisfying $j_1 \neq j_2$, we perturb $F_{j_1}$ and $F_{j_2}$ slightly in a so-called {\it generic} way and investigate the so-called {\it intersection}. This is regarded as the union of circles in ${\rm Int}\ W_f$ and closed intervals whose interiors are in ${\rm Int}\ W_f$ and whose boundaries are in $\partial W_f=q_f(S(f))$, which is connected. The dimension of $F_j$ is $3$, that of $W_f$ is $5$ and the relation $3+3-5=1$ holds. The dimension of the intersection is due to this. 

We can see that the {\it intersection} of ${q_f}^{-1}(F_{j_1})$ and ${q_f}^{-1}(F_{j_2})$ is also canonically obtained. This is regarded as the (disjoint) union of the total spaces of trivial bundles over circles in ${\rm Int}\ W_f$ whose fibers are diffeomorphic to $m-5$ and ($m-4$)-dimensional homotopy spheres. The ($m-4$)-dimensional spheres are the preimages of the closed intervals whose interiors are in ${\rm Int}\ W_f$ and whose boundaries are in $\partial W_f=q_f(S(f))$, which is connected. The ($m-4$)-dimensional spheres are also regarded as the manifolds of the domains of Morse functions with exactly two singular points.
 
 ${\pi}_1(W_f)$ is shown to be trivial in the beginning. Here the homology classes represented by the total spaces of trivial bundles over circles in ${\rm Int}\ W_f$ are zero since $W_f$ is simply-connected. 
The connectedness of $S(f)$ and $q_f(S(f))$ implies the existence of a suitable homotopy $H:S^{m-4} \times [0,1] \rightarrow W_f$. 
This satisfies the following properties. We regard the homotopy spheres of the domains as standard spheres here and this does not affect our arguments.
\begin{itemize}
	\item 
The restriction of the homotopy to $S^{m-4} \times \{0\}$ is regarded as the composition of the embedding into $M$ obtained canonically from each of the ($m-4$)-dimensional homotopy spheres here with $q_f$. 
\item 
The restriction of the homotopy to $S^{m-4} \times \{1\}$ is regarded as some constant map.
\item Let ${\rm S}_{H,\partial W_f}(t)$ be the set of all points in $S^{m-4}$ such that the values of $H$ at the pairs of the points and $t \in [0,1]$ are in the set $\partial W_f$. ${\rm S}_{H,\partial W_f}(t_1) \subset {\rm S}_{H,\partial W_f}(t_2)$ for $0 \leq t_1<t_2 \leq 1$. 
\end{itemize}
We also apply the fact that $W_f$ is simply-connected here or in the third property.

 This has shown the vanishing of the ($m-4$)-th integral homology classes of $M$ represented by the ($m-4$)-dimensional closed manifolds appearing as connected components of the intersection. 

Thus the cup product of $a_{j_1} {h_{j_1}}^{\ast}=a_{j_1} {h_{j_1}}^{\rm c}$ and $a_{j_2} {h_{j_2}}^{\ast}=a_{j_2} {h_{j_2}}^{\rm c}$ vanishes.

We also perturb $F_j$ slightly in a so-called {\it generic} way and investigate the so-called {\it self-intersection} of ${F_j}$. As before, this is also regarded as the union of circles in ${\rm Int}\ W_f$ and closed intervals whose interiors are in ${\rm Int}\ W_f$ and whose boundaries are in $\partial W_f=q_f(S(f))$, which is connected. As before, the self-intersection of ${q_f}^{-1}(F_j)$ gives the vanishing of the ($m-4$)-th integral homology classes of $M$ represented by the ($m-4$)-dimensional closed manifolds appearing as connected components of the self-intersection. Thus the square of $a_j {h_j}^{\ast}=a_j {h_j}^{\rm c}$ vanishes.

$H_{m-4}(M;\mathbb{Z})$ and $H^{4}(M;\mathbb{Z})$ are isomorphic by virtue of Poincar\'e duality theorem. This yields the fact that the cup product of ${h_{j_1}}^{\rm c}$ and ${h_{j_2}}^{\rm c}$ for $j_1 \neq j_2$ and the square of ${h_j}^{\rm c}$ vanish. This completes the proof. \\
\ \\
\underline{Case 2} ${h_j}^{\prime} \in H_2(W_f;\mathbb{Z})$ is not divisible by $2$ and $S(f)$ is not necessarily connected (Main Theorem \ref{mthm:1} (\ref{mthm:1.3.2})). \\
We consider intersections and self-intersections as in Case 1. $W_f$ is spin since it is smoothly immersed into ${\mathbb{R}}^n$ by a smooth immersion of codimension $0$. The embedded submanifold $F_j$ is spin since it is $3$-dimensional, compact and orientable. A normal bundle of $F_j \subset W_f$ is also spin as a result. 

Each of the (self)-intersections of the ($m-2$)-dimensional closed submanifolds with no boundaries in $M$ is, as before, regarded as the union of the total spaces of trivial bundles over circles in ${\rm Int}\ W_f$ whose fibers are diffeomorphic to $S^{m-5}$ and the preimages of closed intervals whose interiors are in ${\rm Int}\ W_f$ and whose boundaries are in $\partial W_f=q_f(S(f))$.

Since a normal bundle of $F_j \subset W_f$ is spin, by a suitable deformation, we can obtain the self-intersection represented as the disjoint union of the following two. 

\begin{itemize}
	\item Smoothly and disjointly embedded circles in ${\rm Int}\ F_j$ satisfying the following two.
	\begin{itemize}
	 \item They are also in ${\rm Int}\ W_f$ and are null-homotopic where the circles are seen as maps into ${\rm Int}\ W_f$ naturally.
	 \item We have the vanishing of the homology class represented by the total spaces of trivial bundles over circles in ${\rm Int}\ W_f$ whose fibers are homotopy spheres: we apply the methods as before.
	\end{itemize}
	\item Smoothly and disjointly embedded connected curves in ${\rm Int}\ F_j$ diffeomorphic to closed intervals satisfying the following three.
	\begin{itemize}
		\item Their interiors are in ${\rm Int}\ F_j$ and ${\rm Int}\ W_f$ and boundaries are in $\partial F_j$ and  $\partial W_f$.
	\item For each class $H_1(W_f;\partial W_f;\mathbb{Z})$, this is represented by exactly two of the connected curves or represented by no curves here.
	\item The preimage of each connected curve is a ($m-4$)-dimensional homotopy sphere: the ($m-4$)-dimensional spheres are also regarded as the manifolds of the domains of Morse functions with exactly two singular points as before. 
\end{itemize}
\end{itemize}
Thus the sum of the integral homology classes represented by these ($m-4$)-dimensional homotopy spheres of the preimages of these closed intervals here is twice some ($m-4$)-th integral homology class of $M$.
%Thanks to this, we can also regard that
%the sum of the homology classes represented by the ($m-4$)-dimensional homotopy spheres of the preimages of the closed intervals here is twice a ($m-4$)-th homology class in $M$.

Thus, by an argument similar to one in Case 1, we have that the square of ${\phi}_{2,\mathbb{Z}}({h_{j}}^{\rm c})$ always vanishes. 
%Note that each ${h_j}^{\prime} \in H_2(W_f;\mathbb{Z})$ is not divisible by $2$ and that this is essential here.

This completes the proof. 
 
\end{proof}

%\begin{Rem}
%\label{rem:3}
%STEP 1 can be regarded as another explanation in the proof of the fact that the 2nd integral homology group of the manifold admitting a special generic map into ${\mathbb{R}}^4$ is free in Theorem \ref{thm:2} (\ref{thm:2.3}) or a main theorem of \cite{nishioka}. \cite{nishioka} first shows that for a special generic map $f:M^5 \rightarrow {\mathbb{R}}^4$ on a $5$-dimensional closed and simply-connected manifold $M^5$ into ${\mathbb{R}}^4$, the 2nd integral homology group $H_2(W_f;\mathbb{Z})$ of $W_f$ is free. After that this shows that the 2nd integral homology group $H_2(M^5;\mathbb{Z})$ is isomorphic to the direct sum of the 3rd integral cohomology group $H^3(W_f;\mathbb{Z})$ and $H_2(W_f;\mathbb{Z})$. We can see that $H_{\ast}(M^5;\mathbb{Z})$ is free.

%%For STEP 2, similar submanifolds appear in the preprint \cite{kitazawa2}.
%\end{Rem}
We can easily see the following by a fundamental topological property of the complex projective space.
\begin{Cor}[Main Corollary \ref{mcor:1}]
	\label{cor:2}
The $3$-dimensional complex projective space does not have any special generic map as in Main Theorem \ref{mthm:1}. If it admits a special generic map into ${\mathbb{R}}^5$, then it does not have any family of homology classes as in Main Theorem \ref{mthm:1}.
\end{Cor}

\begin{Rem}
\label{rem:3}
The $1$-dimensional complex projective plane is regarded as the $2$-dimensional unit sphere. The $2$-dimensional complex projective space does not admit special generic maps into ${\mathbb{R}}^n$ for $n=1,2,3$ by Theorem \ref{thm:2} (\ref{thm:2.1}) and (\ref{thm:2.2}) or Theorem \ref{thm:3}. It does not admit ones into ${\mathbb{R}}^4$ by virtue of the theory \cite{eliashberg} or \cite{kikuchisaeki}. 
Corollary \ref{cor:1} and the theory of \cite{eliashberg} imply that the $3$-dimensional complex projective space does not admit a special generic map into ${\mathbb{R}}^n$ for $n=1,2,3,4,6$. 
%We do not know whether the $3$-dimensional complex projective space admits a special generic map into ${\mathbb{R}}^5$.
\end{Rem}
\begin{Rem}
	\label{rem:4}
	In Main Theorem \ref{mthm:1}, for $m \geq 7$, each integral homology class in $H_j(M;\mathbb{Z})$ and $H_j(W_f;\mathbb{Z})$ satisfying $j \leq m-5$ is represented by a closed, connected and orientable manifold in $M$ and ${\rm Int}\ W_f$. 
	This is due to fundamental arguments on the singularity theory of smooth maps and differential topology.
	More precisely, it is represented by a $j$-dimensional standard sphere. This is due to the assumption that $M$ is simply-connected, the fact that $W_f$ is simply-connected and Hurewicz theorem. \cite{saekisuzuoka} such as Corollary 4.8 and articles \cite{kitazawa0.1,kitazawa0.2,kitazawa0.3} also show that $q_f$ induces the isomorphisms ${q_f}_{\ast}:H_j(M;\mathbb{Z}) \rightarrow H_j(W_f;\mathbb{Z})$, ${q_f}_{\ast}:{\pi}_j(M) \rightarrow {\pi}_j(W_f)$, and ${q_f}^{\ast}:H^j(W_f;\mathbb{Z}) \rightarrow H^j(M;\mathbb{Z})$. These facts on these isomorphisms are more general propositions. 
	We can show them using arguments on handles in our proof of Theorem \ref{thm:3}. They are exercises to readers or we can check related expositions in the articles introduced here.
	Note also that in our case, the $j$-dimensional standard spheres in ${\rm Int}\ W_f$ can be regarded as the images of the $j$-dimensional standard spheres in $M$ where the spheres are ones represented by the homology classes in the beginning.
	
	Throughout the present paper, it seems that we can weaken several conditions. For example, it seems to have no problems to drop the condition that the manifold $M$ and some spaces around it are simply-connected in most cases. Such problems are left to readers, the author and us all. We concentrate on simply-connected manifolds and spaces, mainly closed and simply-connected manifolds whose dimensions are greater than $5$, which are fundamental and important objects in classical and sophisticated algebraic topology and differential topology, in the present paper.
\end{Rem}
Remark \ref{rem:4} with Theorem \ref{thm:3} yields our Main Corollary \ref{mcor:2}.
\begin{Cor}[Main Corollary \ref{mcor:2}]
\label{cor:3}
Let $m \geq 7$ be an integer.
Let $M$ be an $m$-dimensional closed and simply-connected manifold. If there exists a 2nd integral cohomology class whose square is not divisible by $2$, then $M$ admits no special generic maps into $ {\mathbb{R}}^n$ for $n=1,2,3,4,5$.
\end{Cor}

We show some related examples.
\begin{Ex}
\label{Ex}
\label{ex:2}
Let $m \geq 6$ be an integer. Consider a special generic map $f:M \rightarrow {\mathbb{R}}^5$ as in Example \ref{ex:1} on a manifold $M$ represented as a connected sum of $l_0>0$ copies of $S^2 \times S^{m-2}$ whose image is represented as a boundary connected sum of $l_0>0$ copies of $S^2 \times D^3$. This is for the case where $l=l_0$ in Main Theorem \ref{mthm:1} (\ref{mthm:1.3.1}) and (\ref{mthm:1.3.2}).
\end{Ex}
\begin{Ex}
\label{ex:3}
$M:=S^2 \times S^2 \times S^{m-4}$ admits a special generic map $f:M \rightarrow {\mathbb{R}}^5$ satisfying the following two conditions where $m \geq 5$.
\begin{itemize}
\item $f {\mid}_{S(f)}$ is an embedding.
\item $f(M)$ is diffeomorphic to $S^2 \times S^2 \times D^1$.
\end{itemize}
We construct such a map as follows. First we construct the product map of a canonical projection of $S^{m-4}$ to ${\mathbb{R}}={\mathbb{R}}^1$ and the identity map on $S^2 \times S^2$. We embed the manifold of the target suitably into
${\mathbb{R}}^5$. 

Let $m>6$. This is for the case where $l=2$ in Main Theorem \ref{mthm:1} (\ref{mthm:1.3.2}).

Let $m \geq 5$ again. $M$ also admits a special generic map $f_{n-5}:M \rightarrow {\mathbb{R}}^n$ satisfying the following two conditions for any $5 \leq n \leq m$ where $f_0=f$. We have this map in a similar way where we replace the canonical projection by another canonical projection to ${\mathbb{R}}^{n-4}$.
\begin{itemize}
\item $f_{n} {\mid}_{S(f_n)}$ is an embedding.
\item $f_{n}(M)$ is diffeomorphic to $S^2 \times S^2 \times D^{n-4}$.
\end{itemize}
Note that $M$ does not admit special generic maps into ${\mathbb{R}}^n$ for $m \geq 6$ and $n=1,2,3,4$ by Theorem \ref{thm:3}.
\end{Ex}
\begin{Ex}
	\label{ex:4}
	\cite{kreck} recently obtained a classification of $7$-dimensional closed and simply-connected manifolds whose 2nd integral homology groups are free via concrete algebraic topological tools or concrete bordism theory.
	
	\cite{wang} studies classifications of these manifolds which are spin and whose integral cohomology rings are isomorphic to that of the product of the $2$-dimensional complex projective space and $S^3$. \cite{kitazawa0, kitazawa5} construct fold maps on such manifolds. Such manifolds do not admit special generic maps into ${\mathbb{R}}^n$ for $n=1,2,3,4,5$ by virtue of Corollary \ref{cor:3}.
	
	This new discovery is closely related to a question by Takahiro Yamamoto in a contributed talk by the author in  Autumn Meeting of The Mathematical Society of Japan in 2021 (https://www.mathsoc.jp/en/meeting/chiba21sept/). The author introduced \cite{kitazawa0.4}, one of the main theorems of which shows the non-existence of special generic maps into ${\mathbb{R}}^n$ for $n=1,2,3,4,5$ on $7$-dimensional closed and simply-connected manifolds having non-vanishing {\it triple Massey products} (see \cite{kraines, taylor, taylor2} for fundamental theory of {\it triple Massey products}). After the talk, he contacted and asked the author about the inverse of the result. Some observations such as Corollary \ref{cor:3} in the present paper give counterexamples.
\end{Ex}

For some recent examples of special generic maps on closed and simply-connected manifolds, see also \cite{kitazawa2, kitazawa3} for example.

\begin{proof}[A proof of Main Theorem \ref{mthm:2}]
First, in the present situation, we essentially show the case $m=6$ only. The case $m \geq 7$ is not so difficult. This is due to Remark \ref{rem:4}.

$h_0$ is represented by a closed smooth manifold and we may regard this to be diffeomorphic to $S^2$ by Hurewicz theorem: the 2nd homotopy group and the 2nd integral homology group are isomorphic and this isomorphism is given by a canonical isomorphism. Note that $M$ is simply-connected and that this is essential here.

We consider the composition of a suitable smooth embedding $e:S^2 \rightarrow M$ giving the submanifold $e(S^2) \subset M$ with $q_f$.
Suppose that the image of this is a subset of ${\rm Int}\ W_f$. This is an essential assumption in the case $m=6$.
 Suppose also that a 2nd integral homology class in $H_2(W_f;\mathbb{Z})$ of a finite order is realized as the value of 
the homomorphism $(q_f \circ e)_{\ast}$ at the fundamental class where the sphere of the domain is suitably oriented. From this, we have a non-zero integer $a_0$ satisfying the following two.
\begin{itemize} 
\item $a_0\ {q_f}_{\ast}(h_0) \in H_2(W_f;\mathbb{Z})$ is represented by the submanifold given by a smooth embedding of $S^2$ into ${\rm Int}\ W_f$. For example we can consider a suitable connected sum of finitely many embeddings smoothly isotopic to the original smooth embedding $e$ to obtain a desired embedding as the composition with $q_f$.
\item $a_0\ {q_f}_{\ast}(h_0) \in H_2(W_f;\mathbb{Z})$ is zero.
\end{itemize}
$W_f$ is simply-connected as in Main Theorem \ref{mthm:1}. The 2nd property implies that $a_0\ h_0$ also vanishes. This is a contradiction.

This completes the proof.
\end{proof}
We omit examples for Main Theorem \ref{mthm:2}. It seems to discuss examples in systematic studies of closed and simply-connected manifolds whose dimensions are greater than $5$ and whose 2nd integral homology groups are not free. For such manifolds, refer to classification results which will be presented in Remark \ref{rem:5} and some manifolds in \cite{kitazawa2,kitazawa3}.

\begin{proof}[A proof of Main Theorem \ref{mthm:3}]
 We have a situation in Main Theorem \ref{mthm:1} (\ref{mthm:1.3.1}) due to arguments on homology groups and cohomology groups in Remark \ref{rem:4}. This completes the proof.

\end{proof}
\begin{Ex}
The special generic map in Example \ref{ex:2} is for Main Theorem \ref{mthm:3} where the dimension $m>6$ is assumed. $M:=S^2 \times S^2 \times S^{m-2}$ in Example \ref{ex:3} does not admit special generic maps in Main Theorem \ref{mthm:3}.
\end{Ex}

Theorem \ref{thm:3} implies that for $n=1,2,3,4$ the existence of a special generic map $f:M \rightarrow {\mathbb{R}}^n$ on an $m$-dimensional closed manifold $M$ yields the vanishing of the cup product $c_1 \cup c_2$ for any pair of cohomology classes $c_1,c_2 \in H^2(M;\mathbb{Z})$ where $m \geq 6$. Main Theorem \ref{mthm:3} may be regarded as a theorem similar to Theorem \ref{thm:3}. However, it is essentially different from Theorem \ref{thm:3}.

\begin{Rem}
\label{rem:5}
For $6$-dimensional closed and simply-connected manifolds and their classifications, see \cite{jupp,wall,zhubr} for example. See also a homepage \cite{zhubr2}.

The $3$-dimensional complex projective space and the manifold $M$ in Example \ref{ex:3} are also so-called {\it generalized Bott manifolds}. For {\it generalized Bott manifolds}, see \cite{choimasudasuh} for example. They are also so-called {\it toric} manifolds. (Compact) toric manifolds are simply-connected and their diffeomorphism types are conjectured to be determined by their integral cohomology rings.
\end{Rem}
\begin{Rem}
\label{rem:6}
Related to Remarks \ref{rem:3} and \ref{rem:5}, the existence or non-existence of special generic maps on $3$-dimensional complex projective space into ${\mathbb{R}}^5$ is a new open problem.
\end{Rem}

\section{Acknowledgement.}
The author is a member of the project supported by JSPS KAKENHI Grant Number JP17H06128 "Innovative research of geometric topology and singularities of differentiable mappings"
(Principal investigator: Osamu Saeki). The present study is also supported by the project. 
The author would like to thank Takahiro Yamamoto for very natural and interesting questions and his interests in the present study and related studies. 
We declare that data essentially supporting the present study are all in the present article.


\begin{thebibliography}{30}
%\bibitem{akhmetev} P. M. Akhmet'ev, \textsl{On an isotopic and a discrete realization of mappings of an $n$-dimensional sphere in Euclidean space}, Mat. Sb. 187 (1996), 3--34.
\bibitem{barden} D. Barden, \textsl{Simply Connected Five-Manifolds}, Ann. of Math. (3) 82 (1965), 365--385.
\bibitem{burletderham} O. Burlet and G. de Rham, \textsl{Sur certaines applications g\'en\'eriques d'une vari\'et\'e close a $3$ dimensions dans le plan}, Enseign. Math. 20 (1974). 275--292.
\bibitem{calabi} E. Calabi, Quasi-surjective mappings and a generation of Morse theory, Proc. U.S.-Japan Seminar in Differential Geometry, Kyoto, 1965, pp. 13--16.
\bibitem{choimasudasuh} S. Choi, M. Masuda and D. Y. Suh, \textsl{Topological classification of generalized Bott towers}, Trans. Amer. Math. Soc. 362 (2010), 1097--1112.
% \bibitem{cerf} J. Cerf, \textsl{La stratification naturelle des espaces de fonctions deff\'erentiables r'eelles et le th'eor`eme de la pseudo-isotopie}, Inst. Hautes Etudes Sci. Publ. Math. 39 (1970), 5--173.
%\bibitem{crowleyescher} D. Crowley and C. Escher, \textsl{A classification of $S^3$-bundles over $S^4$}, Differential. Geom. Appl. 18 (2003), 363--380, arxiv:0004147.
%\bibitem{crowleynordstrom} D. Crowley and J. Nordstr\"{o}m, \textsl{The classification of $2$-connected $7$-manifolds}, Proc. London. Math. Soc. 119 (2019), 1--54, arxiv:1406.2226.
%\bibitem{dranishnikovrudyak} A. N. Dranishnikov and Y. B. Rudyak, \textsl{Examples of non-formal closed simply connected manifolds of dim ensions 7 and more}, arXiv:math/0306299v3.
%\bibitem{dranishnikovrudyak2} A. N. Dranishnikov and Y. B. Rudyak, \textsl{Examples of non-formal closed {\rm (}$k-1${\rm )}-connected manifolds of dimensions $4k-1$ and more}, arXiv:math/0306299.
\bibitem{eellskuiper} J. J. Eells and N. H. Kuiper, \textsl{An invariant for certain smooth manifolds}, Ann. Mat. Pura Appl. 60 (1962), 93--110.
%\bibitem{ehresmann} C. Ehresmann, \textsl{Les connexions infinitesimales dans un espace fibre differentiable}, Colloque de Topologie, Bruxelles (1950), 29--55.
\bibitem{eliashberg} Y. Eliashberg, \textsl{On singularities of folding type}, Math. USSR Izv. 4 (1970). 1119--1134.
%\bibitem{eliashberg2} Y. Eliashberg, \textsl{Surgery of singularities of smooth mappings}, Math. USSR Izv. 6 (1972). 1302--1326.
%\bibitem{fernandezMunoz} M. Fern\'andez and V. mu\~noz, \textsl{On non-formal simply connected manifolds}, Topology Appl. 135 Issues 1--3 (2004), 111--117. math.DG/0212141.
\bibitem{furuyaporto} Y. K. S. Furuya and P. Porto, \textsl{On special generic maps from a closed manifold into the plane}, Topology Appl. 35 (1990), 41--52.
\bibitem{golubitskyguillemin} M. Golubitsky and V. Guillemin, \textsl{Stable mappings and their singularities}, Graduate Texts in Mathematics (14), Springer-Verlag (1974).
%\bibitem{grossbergkarshon} M. Grossberg and Y. Karshon, \textsl{Bott towers, complete integrability, and the extended character of representations}, Duke Math. J 76 (1994), 23--58.
% \bibitem{haefliger} A. Haefliger, \textsl{Knotted {\rm (}$4k-1${\rm )}-spheres in $6k$-space}, Ann. of Math. 75 (1962), 452--466.
%\bibitem{hatcher} A. E. Hatcher, \textsl{A proof of the Smale conjecture}, Ann. of Math. 117 (1983), 553--607.
\bibitem{hatcher} A. Hatcher, \textsl{Algebraic Topology}, A modern, geometrically flavored introduction to algebraic topology, Cambridge: Cambridge University Press (2002).  
% \bibitem{hiratuka} J. T. Hiratuka, \textsl{A fatorizacao de Stein e o numero de singularidades de aplicacoes estaveis} (in
% Portuguese), PhD Thesis, University of Sao Paulo (2001).
\bibitem{jupp}  P. E. Jupp, \textsl{Classification of certain $6$-manifolds}, Proc. Cambridge Philos. Soc. 73 (1973), 293--300.
%\bibitem{kasuya} H. Kasuya, \textsl{Cohomologically symplectic solvmanifolds are symplectic}, Journal .of Symplectic Geometry 9 (4) (2011), 429--434, arXiv:1005.1157.
\bibitem{kervairemilnor} M. Kervaire and J. W. Milnor, \textsl{Groups of homotopy spheres : I}, Ann. of Math.,, 77 (1963), 504--537.
\bibitem{kikuchisaeki} S. Kikuchi and O. Saeki, \textsl{Remarks on the topology of folds}, Proc. Amer. Math, Soc. No.3 123 (1995), 905--908.
% \bibitem{kitazawa} N. Kitazawa, \textsl{Gluings and decompositions of spaces and maps}, in preparation.
\bibitem{kitazawa0.1} N. Kitazawa, \textsl{On round fold maps} (in Japanese), RIMS Kokyuroku Bessatsu B38 (2013), 45--59.
\bibitem{kitazawa0.2} N. Kitazawa, \textsl{On manifolds admitting fold maps with singular value sets of concentric spheres}, Doctoral Dissertation, Tokyo Institute of Technology (2014).
\bibitem{kitazawa0.3} N. Kitazawa, \textsl{Fold maps with singular value sets of concentric spheres}, Hokkaido Mathematical Journal Vol.43, No.3 (2014), 327--359.
%\bibitem{kitazawa3-1} N. Kitazawa, \textsl{Constructions of round fold maps on smooth bundles}, Tokyo J. of Math. Volume 37, Number 2, 385--403, arxiv:1305.1708.
%\bibitem{kitazawa4} N. Kitazawa, \textsl{Round fold maps and the topologies and the differentiable structures of manifolds admitting explicit ones}, submitted to a refereed journal, arXiv:1304.0618 (the title has changed).%
%\bibitem{kitazawa5} N. Kitazawa, \textsl{Constructing fold maps by surgery operations and homological information of their Reeb spaces}, submitted to a refereed journal, arxiv:1508.05630 (the title has been changed).%
%\bibitem{kitazawa4} N. Kitazawa, \textsl{Notes on fold maps obtained by surgery operations and algebraic information of their Reeb spaces}, submitted to a refereed journal, arxiv:1811.04080.
%\bibitem{kitazawa5} N. Kitazawa, \textsl{Lifts of spherical Morse functions}, submitted to a refereed journal, arxiv:1805.05852.
%\bibitem{kitazawa6} N. Kitazawa, \textsl{Generalizations of Reeb spaces of special generic maps and applications to a problem of lifts of smooth maps}, arxiv:1805.07783.  
%\bibitem{kitazawa7} N. Kitazawa, \textsl{A new explicit way of obtaining special generic maps into the $3$-dimensional Euclidean space}, arxiv:1806.04581.
\bibitem{kitazawa0.4} N. Kitazawa, \textsl{Special generic maps and fold maps and information on triple Massey products of higher dimensional differentiable manifolds}, submitted to a refereed journal, arxiv:2006.08960.
\bibitem{kitazawa0} N. Kitazawa, \textsl{$7$-dimensional simply-connected spin manifolds whose integral cohomology rings are isomorphic to that of ${{\mathbb{C}}P}^2 \times S^3$ admit round fold maps}, submitted to a refereed journal, arxiv:2007.03474v8.
\bibitem{kitazawa1} N. Kitazawa, \textsl{Closed manifolds admitting no special generic maps whose codimensions are negative and their cohomology rings}, submitted to a refereed journal, arxiv:2008.04226v5.
\bibitem{kitazawa2} N. Kitazawa, \textsl{Notes on explicit special generic maps into Euclidean spaces whose dimensions are greater than $4$}, a revised version will be submitted based on positive comments (major revision) by referees and editors after the first submission to a refereed journal, arxiv:2010.10078. 
\bibitem{kitazawa3} N. Kitazawa, \textsl{The images of special generic maps of several classes}, arxiv:2011.12066.
\bibitem{kitazawa4} N. Kitazawa, \textsl{$7$-dimensional closed simply-connected and spin manifolds having 2nd integral cohomology classes whose squares are not divisible by $2$ and stable fold maps on them}, arXiv:2104.10871.
\bibitem{kitazawa5} N. Kitazawa, \textsl{A note on cohomological structures of special generic maps}, a revised version is submitted based on positive comments by referees and editors after the second submission to a refereed journal.
%\bibitem{kitazawasaeki} N. Kitazawa and O. Saeki, \textsl{Round fold maps of $n$-dimensional manifolds into ${\mathbb{R}}^{n-1}$}, arXiv:2111.13510.
%\bibitem{kitazawa8} N. Kitazawa \textsl{Surgery operations to fold maps to construct fold maps whose singular value sets may have crossings}, arxiv:2003.04147.
%\bibitem{kitazawa9} N. Kitazawa \textsl{Surgery operations to fold maps to increase connected components of singular sets by two}, arxiv:2004.03583.
% \bibitem{kobayashi} M. Kobayashi, \textsl{Simplifying certain stable mappings from simply connected %R% $4$-manifolds into the plane}, Tokyo J. Math. 15 (1992), 327--349.
%\bibitem{kobayashi} M. Kobayashi, \textsl{Stable mappings with trivial monodromies and application to inactive log-transformations}, RIMS Kokyuroku. 815 (1992), 47--53.
%\bibitem{kobayashi2} M. Kobayashi, \textsl{Bubbling surgery on a smooth map}, preprint.
%\bibitem{kobayashisaeki} M. Kobayashi and O. Saeki, \textsl{Simplifying stable mappings into the plane from a global viewpoint}, Trans. Amer. Math. Soc. 348 (1996), 2607--2636. 
%
\bibitem{kraines} D. Kraines, \textsl{Massey higher products}, Trans. Amer. Math. Soc. 124 (1966), 431--449. 
\bibitem{kreck} M. Kreck,  \textsl{On the classification of $1$-connected $7$-manifolds with torsion free second homology}, to appear in the Journal of Topology, arxiv:1805.02391.
% \bibitem{kurokisuh} S. Kuroki and D. Y. Suh, \textsl{Cohomological non-rigidity of eight-dimensional complex projective
% towers}, to appear in Algebraic \& Geometric Topology; OCAMI preprint series 13-12.
%\bibitem{masseyuehara} W. Massey and H. Uehara, \textsl{the Jacobi identity for Whitehead products}, Algebraic geometry and topology. A symposium in honor of S. Lefschetz, pp. 361--377, Princeton University Press, Princeton, N. J., 1957.
%\bibitem{luptonoprea} G. Lupton and J. Oprea, \textsl{Cohomologically Symplectic spaces: toral actions and the Gotilieb group}, Trans. Amer. Math. Soc. 347 (1) (1995), 261--288. 
%\bibitem{masuda} M. Masuda, \textsl{Classification of real Bott manifolds}, arxiv:0809.2178.
\bibitem{milnor} J. Milnor, \textsl{On manifolds homeomorphic to the $7$-sphere}, Ann. of Math. (2) 64 (1956), 399--405.
\bibitem{milnor2} J. Milnor, \textsl{Morse Theory}, Annals of Mathematic Studies AM-51, Princeton University Press; 1st Edition (1963.5.1).
\bibitem{milnor3} J. Milnor, \textsl{Lectures on the h-cobordism theorem}, Math. Notes, Princeton Univ. Press, Princeton, N.J. 1965.
\bibitem{milnorstasheff} J. Milnor and J. Stasheff, \textsl{Characteristic classes}, Annals of Mathematics Studies, No. 76. Princeton, N. J; Princeton University Press (1974).
%\bibitem{murai} T. Murai, \textsl{{\rm (}Master's thesis : in Japanese{\rm )}}, Master's thesis, Tsuda. Univ.
%\bibitem{nakamura} K. Nakamura, K. Nakamura, \textsl{{\rm (}Bachelor's thesis : in Japanese{\rm )}}, Hirosdhima Univ. 
\bibitem{nishioka} M. Nishioka, \textsl{Special generic maps of $5$-dimensional manifolds}, Revue Roumaine de Math`{e}matiques Pures et Appliqu\`{e}es, Volume LX No.4 (2015), 507--517.
%\bibitem{reeb} G. Reeb, \textsl{Sur les points singuliers d\`{u}ne forme de Pfaff completement integrable ou d’une fonction numerique}, -C. R. A. S. Paris 222 (1946), 847--849. 
% \bibitem{saeki} O. Saeki, \textsl{Notes on the topology of folds}, J. Math. Soc. Japan Volume 44, Number 3 (1992), 551--566.
\bibitem{saeki} O. Saeki, \textsl{Topology of special generic maps of manifolds into Euclidean spaces}, Topology Appl. 49 (1993), 265--293.
\bibitem{saeki2} O. Saeki, \textsl{Topology of special generic maps into $\mathbb{R}^3$}, Workshop on Real and Complex Singularities (Sao Carlos, 1992), Mat. Contemp. 5 (1993), 161--186.
% \bibitem{saeki4} O. Saeki, \textsl{Morse functions with sphere fibers}, Hiroshima Math. J. Volume 36, Number 1 (2006), 141--170.
\bibitem{saekisakuma} O. Saeki and K. Sakuma, \textsl{On special generic maps into ${\mathbb{R}}^3$}, Pacific J. Math. 184 (1998), 175--193.
\bibitem{saekisakuma2} O. Saeki and K. Sakuma, \textsl{Special generic maps of $4$-manifolds and compact complex analytic surfaces}, Math. Ann. 313, 617--633, 1999.
\bibitem{saekisuzuoka} O. Saeki and K. Suzuoka, \textsl{Generic smooth maps with sphere fibers} J. Math. Soc. Japan Volume 57, Number 3 (2005), 881--902.
\bibitem{sakuma} K. Sakuma, \textsl{On special generic maps of simply connected $2n$-manifolds into ${\mathbb{R}}^3$},Topology Appl. 50 (1993), 249--261.
%\bibitem{sakuma2} K. Sakuma, \textsl{On the topology of simple fold maps}, Tokyo J. of Math. Volume 17, Number 1 (1994), 21--32.
%\bibitem{shiota} M. Shiota, \textsl{Thom's conjecture on triangulations of maps}, Topology 39 (2000), 383--399.
%\bibitem{smale}, S. Smale, \textsl{Diffeomorphisms of the 2-sphere}, Proc. Amer. Math, Soc. 10 (1959), 621--626.
% \bibitem{smale} S. Smale, \textsl{Generalized Poincare's conjecture in dimensions greater than four}, Ann. of Math. (2) 74 (1961) 391--406.
% \bibitem{smale2} S. Smale, \textsl{On the structure of manifolds}, Amer. J. Math 84 (1962), 387--399.
\bibitem{steenrod} N. Steenrod, \textsl{The topology of fibre bundles}, Princeton University Press (1951). 
%\bibitem{suzuoka} K. Suzuoka, \textsl{On the topology of sphere type fold maps of $4$-manifolds into the plane}, PhD Thesis, Univ. of Tokyo (March 2002).
\bibitem{taylor} L. R. Taylor, \textsl{Controlling indeterminacy in Massey triple products}, Geom. Dedicata 148 (2010), 371--389. 
\bibitem{taylor2} L. R. Taylor, \textsl{Massey Triple Products}, https://www3.nd.edu/\~{}taylor/talks/2011-03-22-Princeton.pdf, Princeton Topology Seminar, 2011/3/22.
\bibitem{thom} R. Thom, \textsl{Quelques propri\'et\'es globales des vari\'et\'es diff\'erentiables}, Commentarii Mathematici Helvetici (1954), Volume 28, 17--86.
%\bibitem{thom} R. Thom, \textsl{Les singularites des applications differentiables}, Ann. Inst. Fourier (Grenoble) 6 (1955-56), 43--87.
%\bibitem{wall} C.T.C. Wall, \textsl{All $3$-manifolds embed in $5$-space}, Bull. Amer. Math. Soc. 71 (1965), 564--567.
% \bibitem{wall} C. T. C. Wall, \textsl{Classification of {\rm (}$n-1${\rm )}-Connected $2n$-Manifolds}, Ann. of Math. Second Series (1) 75 (Jan.1962), 163--189.
\bibitem{wall} C. T. C. Wall, \textsl{Classification problems in differential topology. V. On certain $6$-manifolds}, Invent. Math. 1 (1966), 355--374. 
\bibitem{wang} X. Wang \textsl{On the classification of certain $1$-connected $7$-manifolds and related problems}, arXiv:1810.08474.
%\bibitem{whitney} H. Whitney, \textsl{On singularities of mappings of Euclidean spaces: I, mappings of the plane into the plane}, Ann. of Math. 62 (1955), 374--410.
\bibitem{wrazidlo} D. J. Wrazidlo, \textsl{Standard special generic maps of homotopy spheres into Eucidean spaces}, Topology Appl. 234 (2018), 348--358, arxiv:1707.08646.
\bibitem{wrazidlo2} D. J. Wrazidlo (including a joint work with O. Saeki and K. Sakuma), \textsl{The Milnor $7$-sphere does not admit a special generic map into ${\mathbb{R}}^3$},  http://math.ac.vn/conference/FJV2018/images/slides/Wardzilo\_FJV2018.pdf, 2018. 
\bibitem{wrazidlo3} D. J. Wrazidlo, \textsl{On special generic maps of rational homology spheres into
 Euclidean spaces}, arxiv:2009.05928.
%\bibitem{yamamoto} M. Yamamoto, \textsl{Lifting a generic map of a surface into the plane to an embedding into $4$-space}, Illinois Journal of Mathematics 51 (2007), 705--721.
\bibitem{zhubr} A. V. Zhubr, Closed simply-connected six-dimensional manifolds: proofs of classification theorems, Algebra i Analiz 12 (2000), no. 4, 126--230.
\bibitem{zhubr2} A. V. Zhubr (responsible for the page), http://www.map.mpim-bonn.mpg.de/6-manifolds:\_1-connected.
\end{thebibliography}
\end{document}